\documentclass[12pt]{amsart}

\usepackage{graphicx}
\usepackage{latexsym}
\usepackage[centertags]{amsmath}
\usepackage{amsfonts}
\usepackage{amssymb}
\usepackage{amsthm}
\usepackage{newlfont}
\usepackage[square]{natbib}

\def\mytopsep{3mm}

\newtheoremstyle{myplain}{\mytopsep}{\mytopsep}{\itshape}{0pt}{\bfseries}{.}{3mm}{}
\newtheoremstyle{mydefinition}{\mytopsep}{\mytopsep}{\normalfont}{0pt}{\bfseries}{.}{3mm}{}

\newtheoremstyle{myremark}{\mytopsep}{\mytopsep}{\normalfont}{0pt}{\bfseries}{.}{3mm}{}

\theoremstyle{myplain}

\newtheorem{thm}{Theorem}[section]
\newtheorem{cor}[thm]{Corollary}
\newtheorem{lem}[thm]{Lemma}
\newtheorem{prop}[thm]{Proposition}

\theoremstyle{mydefinition}
\newtheorem{dfn}[thm]{Definition}
\theoremstyle{myremark}
\newtheorem{rem}[thm]{Remark}
\newtheorem{exa}[thm]{Example}

\allowdisplaybreaks[1]

\makeatletter\@addtoreset{equation}{section}\makeatother

\def\rank{\mathop{\mathrm{rank}}}

\def\supp{\mathop{\mbox{supp}}}
\def\ord{\mathrm{ord}}
\def\mb{\mathbf}

\def\NN{\mathbb{N}}

\def\QQ{\mathbb{Q}}
\def\PP{\mathbb{P}}
\def\xx{x^{-1}}

\def\CC{\mathbb{C}}
\def\ZZ{\mathbb{Z}}

\def\ct{\mathop{\mathrm{CT}}}
\def\PT{\mathop{\mathrm{PT}}}
\def\NT{\mathop{\mathrm{NT}}}
\def\pt{\mathop{\mathrm{PT}}}
\def\nt{\mathop{\mathrm{NT}}}
\def\CT{\mathop{\mathrm{CT}}}

\def\uprho{\mbox{}^\rho \,}

\renewcommand{\ll}{\langle\!\langle}
\renewcommand{\gg}{\rangle\!\rangle}

\begin{document}

\title[A Generalization of Stanley's Monster Reciprocity Theorem]{A Generalization of Stanley's Monster Reciprocity Theorem}

\author{Guoce Xin}
\address{Department
of Mathematics\\
Brandeis University\\
Waltham MA 02454-9110} \email{guoce.xin@gmail.com}

\date{April 21, 2005}
\begin{abstract}
By studying the reciprocity property of linear Diophantine systems
in light of Malcev-Neumann series, we present in this paper a new
approach to and a generalization of Stanley's monster reciprocity
theorem. A formula for the ``error term" is given in the case when
the system does not have the reciprocity property. We also give a
short proof of Stanley's reciprocity theorem for linear
homogeneous Diophantine systems.
\end{abstract} \maketitle


\section{Introduction}
Let $A$ be an $r$ by $n$ matrix with integer entries, and let
$\mb{b}$ be an $r$-vector in $\ZZ^r$. Many combinatorial problems
turn out to be equivalent to finding all nonnegative integral
(column) vectors $\alpha\in \NN^n$ satisfying
\begin{equation}\label{e-lid} A\alpha =\mb{b},
\end{equation}
especially in the homogeneous case when $\mb{b}$ equals $\mb{0}$,
of which the solution space is a rational cone. Such problems are
also known as solving a linear Diophantine system.

There are two closely related generating functions associated to
\eqref{e-lid}:
\begin{equation*}
 E(\mb{x};\mb{b})=\sum_{(\alpha_1,\dots
,\alpha_n)\in \NN^n} x_1^{\alpha_1}\cdots x_n^{\alpha_n},\qquad
\bar{E}(\mb{x};\mb{b})=\sum_{(\alpha_1,\dots ,\alpha_n)\in
\PP^n}x_1^{\alpha_1}\cdots x_n^{\alpha_n},
\end{equation*}
where the first sum ranges over all $\alpha=(\alpha_1,\dots
,\alpha_n)$ such that $A\alpha =\mb{b}$, and the second sum ranges
over all positive integral $\alpha$ such that $A\alpha=-\mb{b}$.
We omit $\mb{b}$ in the homogeneous case. The following well-known
reciprocity theorem for homogeneous linear diophantine equations
was given by Stanley as \cite[Theorem 4.1]{stanley-magiclabel}.

\begin{thm}[Reciprocity Theorem] 
\label{t-3-recip}
Let $A$ be an $r$ by $n$ integral matrix of full rank $r$. If
there is at least one $\alpha\in \PP^n$ such that $A\alpha
=\mb{0}$, then we have as rational functions
\begin{equation*}
    E(x_1,\dots,x_n)= (-1)^{n-r} \bar{E}(x_1^{-1},\dots ,x_n^{-1}).
\end{equation*}
\end{thm}

Previous proofs of this theorem used decompositions into
{simplicial cones} or lattice cones, or complicated algebraic
technique. See \cite[p.~214]{stanley-rec} and \cite{stanley-local}
for further information. 
We will give a short proof using a signed cone decomposition and
induction.

In the general situation, the best known result (up to now) is the
monster reciprocity theorem, which was given by Stanley
\cite{stanley-rec} in 1974. The theorem will be stated later after
new notation is introduced. It includes as special cases many
combinatorial reciprocity theorems, such as the reciprocity
theorem for homogeneous linear Diophantine system, that for
Ehrhart polynomials, and that for P-partitions, etc. We will give
a simple approach to this theorem. As applications, we give
detailed, and short, implication of the reciprocal domain theorem
\cite[Proposition~8.3]{stanley-rec}.

The new approach uses the idea of Malcev-Neumann series
\cite{passmann,xinresidue,xinthesis}, which defines a total
ordering on the group of monomials to clarify the series expansion
of rational functions. We study the reciprocity property of
 an object that is more
general, but less combinatorial, than that was studied in
\cite{stanley-rec}. The new objects we are going to study are
\emph{Elliott-rational functions}, while the previous objects are
Elliott-rational functions with a monomial numerator. By an
Elliott-rational function, we mean the one that can be written as
$$ F(\lambda_1,\dots,\lambda_r,\mb{x})= \frac{p(\lambda_1,\dots,\lambda_r,\mb{x})}{
\prod_{i=1}^m (y_i-z_i)},$$ where $p$ is a polynomial and $y_i$
and $z_i$ are monomials.

In this larger set of objects, it is much easier to build up the
reduction steps. Theorem \ref{t-rec-erroterm}, a general result
that gives a reciprocity formula for Elliott-rational functions,
turns out to be easy to prove. We shall use this result to
formulate the monster reciprocity theorem (Theorem
\ref{t-monster}).

In Section 2, we introduce the basic idea of Malcev-Neumann series
and reformulate the reciprocity of linear Diophantine system in
terms of constant terms. In Section 3, we develop the reciprocity
theorem for Elliott-rational functions. We apply our result in
Section 4 to give the generalized monster reciprocity theorem. In
section 5, we illustrate the monster reciprocity theorem by
examples, and as an application, we give a simple derivation of
Theorem \ref{t-3-recip}. Section 6 includes an inductive
(combinatorial) proof of Theorem \ref{t-3-recip}.

\section{Reciprocity in Terms of Constant Terms}

Solving a linear Diophantine system (LD-system for short) means
finding all vectors $\alpha\in \NN^n$ that satisfy
$A\alpha=\mb{b}$, where $A$ is an $r$ by $n$ matrix with integral
entries. More precisely, we want to solve the following system of
equations:
\begin{align}\label{e-ax-b}
a_{1,1}\alpha_1 +a_{1,2} \alpha_2+\cdots +a_{1,n} \alpha_n &=b_1 \nonumber \\
a_{2,1}\alpha_1 +a_{2,2} \alpha_2+\cdots +a_{2,n} \alpha_n &=b_2 \nonumber \\
\cdots \cdots \qquad &= \cdots \\
a_{r,1}\alpha_1 +a_{r,2} \alpha_2+\cdots +a_{r,n} \alpha_n &=b_r.
\nonumber
\end{align}
We assume the rank of $A|\mb{b}$ equals the rank of $A$, for
otherwise, the LD-system has no solution even in $\QQ$.

Let $C_i$ be the $i$th column vector of $A$. Then the above system
is the same as
$$C_1 \alpha_1+C_2\alpha_2+\cdots +C_n\alpha_n=\mb{b}.$$

Now let $E(\mb{b})$ and $\bar{E}(\mb{b})$ be the sets of all such
solutions in $\NN^n$ and $\PP^n$ respectively. It is interesting
to study the following two associated generating functions of
\eqref{e-ax-b}:
\begin{align}
E(\mb{x};\mb{b})&=\sum_{\alpha\in E(\mb{b})} \mb{x}^\alpha,\quad
\qquad
    \bar{E}(\mb{x};\mb{b})= \sum_{\alpha\in \bar{E}(-\mb{b})} \mb{x}^\alpha
\end{align}
where $\mb{x}=(x_1,\dots ,x_n)$ and if $\alpha=(\alpha_1,\dots
,\alpha_n)$, then $\mb{x}^\alpha:= x_1^{\alpha_1}\cdots
x_n^{\alpha_n}.$

The above equation defines two rational functions in $\mb{x}$. If
as rational functions
$E(\mb{x};\mb{b})=(-1)^{n-r}\bar{E}(\mb{x^{-1}};\mb{b})$, then we
say that the system \eqref{e-ax-b} has the $R$-\emph{property}
(short for \emph{reciprocity property}).

We can compute $E(\mb{x};\mb{b})$ by replacing the $r$ linear
constraints with $r$  new variables $\lambda_1, \lambda_2,\dots
,\lambda_r$ and then take the constant terms. Let $\Lambda$ be
$(\lambda_1,\dots ,\lambda_r)$, and let $\ct_\Lambda F$ be the
constant term of $F$ in $\Lambda$. We have
\begin{align}
E(x;\mb{b}) &=  \sum_{\alpha\in \NN^n} \ct_\Lambda
\lambda_1^{a_{1,1} \alpha_1+\cdots +a_{1,n}\alpha_n-b_1}\cdots
\lambda_r^{a_{r,1} \alpha_1+\cdots +a_{r,n}\alpha_n-b_r} \mathbf{x}^\alpha \nonumber\\
&= \ct_\Lambda \frac{\lambda_1^{-b_1}\cdots
\lambda_r^{-b_r}}{\prod_{i=1}^n (1-
\lambda_1^{a_{1,i}}\lambda_2^{a_{2,i}}\cdots \lambda_r^{a_{r,i}}
x_i)}=\ct_\Lambda
\frac{\Lambda^{-\mb{b}}}{\prod_{i=1}^n(1-\Lambda^{C_i}x_i)},
\label{e-mac-Ex}
\end{align}
with the working ring $\CC[\Lambda,\Lambda^{-1}][[\mb{x}]]$, where
$\Lambda^{-1}$ means $(\lambda_1^{-1},\dots ,\lambda_r^{-1})$. The
above conversion can be trait back to MacMahon \cite{mac}.
Similarly we get
\begin{align}
\label{e-mac-Ebx} \bar{E}(\mathbf{x};\mb{b})=\ct_\Lambda
\frac{\Lambda^{\mb{b}}  \prod_{i=1}^n
\Lambda^{C_i}x_i}{\prod_{i=1}^n(1- \Lambda^{C_i}x_i)}.
\end{align}

We define $\mathcal{E}(\Lambda,\mb{x};\mb{b})$ and
$\bar{\mathcal{E}}(\Lambda,\mb{x};\mb{b})$ to be the crude
generating functions of $E(\mathbf{x},\mb{b})$ and
$\bar{E}(\mb{x};\mb{b})$ as
\begin{align}
    \mathcal{E}(\Lambda,\mb{x};\mb{b})
    =
    \frac{\Lambda^{-\mb{b}}}{\prod_{i=1}^n(1-\Lambda^{C_i}x_i)},
    \qquad
    \bar{\mathcal{E}}(\Lambda,\mathbf{x};\mb{b})=
\frac{\Lambda^{\mb{b}}  \prod_{i=1}^n
\Lambda^{C_i}x_i}{\prod_{i=1}^n(1- \Lambda^{C_i}x_i)}
   ,
\end{align}
and observe that as rational functions
$$\bar{\mathcal{E}}(\Lambda^{-1},\mb{x^{-1}};\mb{b})=(-1)^n \mathcal{E}(\Lambda,\mb{x};\mb{b}).$$
However, the series expansion of the two sides of the above
equation are different. The change of variables by $\Lambda\to
\Lambda^{-1}$, which corresponds to multiplying each row of
\eqref{e-ax-b} by $-1$, will not make a difference when taking
constant terms. Therefore, the system has the R-property if and
only if as rational functions
$$ \ct_{\Lambda} \mathcal{E}(\Lambda,\mb{x}) =(-1)^r {\ct_\Lambda}^\prime \mathcal{E}(\Lambda,\mb{x}),$$
where we expand $\mathcal{E}(\Lambda,\mb{x})$ on the LHS at
$\mb{x}=\mb{0}$, while on the RHS at $\mb{x}=\infty$.

As we shall see later, the different expansions appearing in the
above equation is easily explained in the context of
Malcev-Neumann series.

The group of monomials in $\Lambda$ and $\mb{x}$ can be given a
total ordering ``$\preceq^\rho $" that is compatible with its
group structure; i.e., for any monomials $A,B$ and $C$, $A\preceq B$ implies
$AC\preceq^\rho BC$. This is equivalent to a total ordering
$\le^\rho$ on the additive group $\ZZ^{n+r}$. An important such
ordering $\le$ is the reverse lexicographical ordering on
$\ZZ^{n+r}$. Then a \emph{Malcev-Neumann series} (or
\emph{MN-series} for short) with respect to $\preceq^\rho$ is a
formal series on $\Lambda$ and $\mb{x}$ with a well-ordered
\emph{support}: the set of monomials corresponds to the nonzero
terms. Recall that a \emph{well-ordered set} is a totally ordered
set such that every nonempty subset has a minimum.

For our purpose, $\rho$ will denote an injective endomorphism of
$\ZZ^{n+r}$ (a nonsingular integral matrix), and $\le ^\rho$ will
be the induced total ordering defined by $a\le^\rho b$ if and only
if $\rho(a)\le \rho(b)$. We denote by $\CC^\rho\ll
\Lambda,\mb{x}\gg$ the corresponding field of MN-series with
respect to $\rho$. The field of iterated Laurent series $\CC\ll
\Lambda,\mb{x} \gg$, where $\rho$ is the identity map and is
omitted, has been studied in \cite{xiniterate,xinthesis}. For a
more general setting of MN-series, the readers are referred to
\cite{xinresidue,xinthesis} or \cite[Chapter~13]{passmann}.

The series expansion of MN-series will be explained in more
details in the next section. Let us review some properties of
MN-series \cite{xinresidue} to see that such fields are suitable
for dealing with different kinds of series expansions of rational
functions.

For any total ordering $\le^\rho$, $\CC^\rho \ll \Lambda, \mb{x}
\gg$ is a field. In particular, $\CC\ll \Lambda, \mb{x}\gg$ is the
field of iterated Laurent series \cite{xiniterate}.

The field $\CC(\Lambda,\mb{x})$ of rational functions is naturally
embedded into $\CC^\rho \ll \Lambda, \mb{x} \gg$ for any $\rho$.
This follows from the field structure of $\CC^\rho\ll \Lambda,
\mb{x}\gg$ and the fact that every polynomial has a finite
support.

Every rational function $F(\Lambda, \mb{x})$ has a unique
expansion in $\CC^\rho\ll \Lambda, \mb{x} \gg$. The expansions of $F$
for different $\rho$ are usually different. For instance, the
expansion of $1/(x-y)$ in $K\ll x,y\gg $ is
$$\frac{1}{x-y}=\frac{1}{x}\cdot \frac{1}{1-y/x}=\frac{1}{x}\sum_{k\ge 0} y^k/x^k,$$
but the expansion in $K\ll
y,x \gg$ is
$$\frac{1}{x-y}=\frac{1}{-y}\cdot \frac{1}{1-x/y}=\frac{1}{-y}\sum_{k\ge 0} x^k/y^k.$$
Note
that we can write $K\ll y,x\gg$ as $K^\rho \ll x,y\gg$ where
$\rho$ is defined by the matrix $\left(\begin{array}{cc}
  0 & 1 \\
  1 & 0
\end{array} \right)$, or by abuse of notation, $\rho(x)=y$ and
$\rho(y)=x$.

Recall also that every subset of a well-ordered set is
well-ordered. Thus the following operators $\ct_\lambda$,
$\pt_\lambda$, and $\nt_\lambda$ are well-defined for MN-series.
\begin{align*}
\ct_\lambda \sum_{k\in \ZZ} b_k \lambda^k &= b_0,\quad \pt_\lambda
\sum_{k\in \ZZ} b_k \lambda^k = \sum_{k\ge 0} b_k
\lambda^{k},\quad \text{ and } \nt_\lambda \sum_{k\in \ZZ} b_k
\lambda^k = \sum_{k<0} b_k \lambda^{k}.
\end{align*}
Obviously, for an MN-series $F(\lambda)$, $\ct_\lambda F(\lambda)
=\left. \pt_\lambda F(\lambda) \right|_{\lambda=0}$. The constant
term operators are commutative so that taking the constant term in
a set of variables is defined by iteration.

Now it is easy to see that Theorem \ref{t-3-recip} is a
consequence of the following proposition.
\begin{prop}\label{p-3-recip}
Suppose that $\bar{E}$ is nonempty. Then
\begin{align}\label{e-3-recip1}
    \ct_\Lambda  \mathcal{E}(\mathbf{x};\mb{0}) =(-1)^{\mathrm{rank}(A)}\ct_\Lambda \uprho \mathcal{E}(\mathbf{x};\mb{0}),
    \end{align}
    where $\rho$ is the endomorphism defined by $\rho(x_i)=x_i^{-1}$ and
    $\rho(\lambda_i)=\lambda_i$.
\end{prop}
On the other hand, it is easy to deal with the case of $\rank(A)<
r$. So Theorem \ref{t-3-recip} is equivalent to Proposition
\ref{p-3-recip}, whose proof will be given in section 6.

The following lemma asserts that elementary row operation will not
change the solution space of an LD-system. We give it here to show
that all the work can be done algebraically.

\begin{lem}[\cite{xinresidue}, Corollary~3.18]\label{l-3-monomial}
Suppose $\mb{y}$ is another set of variables. If $\Phi \in
K[\mb{x},\mb{\xx}]\ll \mb{y}\gg$, then for $f_i=x_1^{b_{i1}}\cdots
x_n^{b_{in}}$ with $\det(b_{ij})_{1\le i,j\le n}\ne 0$,
$$\ct_{\mb{x}} \Phi(f_1,\dots ,f_n)=\ct_{\mb{x}} \Phi(x_1,\dots ,x_n).$$
\end{lem}

\section{Reciprocity of Elliott-Rational Functions}

It is convenient for our purpose to denote by $K$ the field
$\CC(\mb{x})$. The field of rational functions
$\CC(\Lambda,\mb{x})$ can be identified with $K(\Lambda)$. Usually
we are taking constant terms in the ring
$\CC[\Lambda,\Lambda^{-1}][[\mb{x}]]$. This ring can be embedded
into $\CC^\rho \ll \Lambda, \mb{x} \gg$, as long as $\preceq
^\rho$ is compatible with the relation $x_i \succ\!\!\succ^\rho
\lambda_j$ for all $i$ and $j$, where $A\succ\!\!\succ^{\rho} B$
means that $A\succ^\rho B^k$ for any positive integer $k$.

The case $r=1$ is illustrative for our understanding of the series
expansion for MN-series, and in this particular case, we need not
restrict ourselves to Elliott-rational functions. Let us consider
the following problem.

\vspace{3mm} \noindent \emph{Problem:} Given
 a rational function $Q(\lambda)$ (short for $Q(\lambda,\mb{x})$)
 of $\lambda$ and $\mb{x}$,
compute $\PT^\rho _\lambda Q(\lambda,\mb{x})$, where the notation
$\PT^\rho_\lambda$ indicates that $Q(\lambda)$ is treated as an
element of $\CC^\rho \ll \lambda,\mb{x} \gg$, and we use similar
notations for the $\ct$
 and $\nt$ operators.

To deal with this problem, we shall understand that $Q(\lambda)$
is not only an element of $K(\lambda)$, but also an element of
$\CC^\rho\ll \lambda,\mb{x}\gg$. As an element of $K(\lambda)$,
$Q(\lambda)$ can be written as $p(\lambda)/q(\lambda)$, where
$p(\lambda)$ and $q(\lambda)$ are both in $K[\lambda]$. As an
element of $\CC^\rho\ll \lambda,\mb{x}\gg$, the denominator
$q(\lambda)$ plays an important role.

Recall that $\CC^\rho\ll \lambda,\mb{x}\gg$ is equipped with a
total ordering $\preceq^\rho$ on its group of monomials and that
its elements have well-ordered supports. Thus for a nonzero
element $\eta$, we can define its \emph{order} $\ord\, \eta$ to be
$\min \supp(\eta)$, and its \emph{initial term} to be the term
with the least order. The order of $0$ is treated as $\infty$. Let
us write $q(\lambda)=\sum_{i=0}^{d} a_i \lambda^i$, with $a_i\in
\CC(\mb{x})$ and $a_d\ne 0$. To expand $Q(\lambda)$ into a series
in $\CC^\rho \ll \lambda,\mb{x} \gg$, we need to find the
$\lambda$-{\em initial term} $a_j\lambda^j$, i.e., the $j$ such
that $\ord(a_j \lambda^j)\prec^\rho \ord(a_i\lambda^i)$ for all
$i\ne j$. This can be achieved because of the different powers in
$\lambda$. Then
$$\frac{1}{q(\lambda)} =\frac1{a_j\lambda^j} \frac{1}{1+\sum_{i\ne j}
a_i/a_j\lambda^{i-j} } =\frac1{a_j\lambda^j} \sum_{k\ge 0} (-1)^k
\Big(\sum_{i\ne j} a_i/a_j\lambda^{i-j}\Big)^k.$$ This expansion
is justified by the composition law \cite[Theorem~2.2]{xinthesis}.

It is now clear that we have the following three situations.
\begin{enumerate}
\item If $j$ equals $0$, then for any polynomial $p(\lambda)$,
$p(\lambda)/q(\lambda)$ contains only nonnegative powers in
$\lambda$. In this case, we say that $1/q(\lambda)$ is $\PT^\rho$
in $\lambda$.

\item If $j$ equals $d$, then for any polynomial $p(\lambda)$ of degree
in $\lambda$ less than $d$, $p(\lambda)/q(\lambda)$ contains only
negative powers in $\lambda$. In this case, we say that
$1/q(\lambda)$ is $\NT^\rho$ in $\lambda$.

\item If $j$ equals neither $0$, nor $d$, then
 $1/q(\lambda)$ contains both positive and
negative powers in $\lambda$. Thus $1/q(\lambda)$ is neither
$\PT^\rho$ nor  $\NT^\rho$ in $\lambda$.
\end{enumerate}

\begin{lem}\label{l-3-property}
Let $q_1$ and $q_2$ be polynomials in $\lambda$. Then for any
total ordering $\le^\rho$
\begin{itemize}
    \item Both $1/q_1(\lambda)$ and
$1/q_2(\lambda)$ are $\pt^\rho$ in $\lambda$ if and only if $1/(q
_1q_2)$ is.
    \item Both $1/q_1(\lambda)$ and
$1/q_2(\lambda)$ are $\nt^\rho$ in $\lambda$ if and only if $1/(q
_1q_2)$ is.
    \item For all the other cases, $1/(q_1q_2)$ is neither $\pt^\rho $
    in $\lambda$ nor $\nt^\rho$ in $\lambda$.
\end{itemize}
\end{lem}
\begin{proof}
We prove the first case for $\pt$ as follows. The other cases are
similar. Write $$q_1=\sum_{i=0}^{d_1} a_i \lambda^i, \quad
q_2=\sum_{i=0}^{d_2} b_i \lambda^i, \quad \text{ and
}q_1q_2=\sum_{i=0}^{d_1+d_2} c_i \lambda^i.$$ Suppose that
$a_{j_1}\lambda^{j_1}$ and $b_{j_2}\lambda^{j_2}$ are the
$\lambda$-initial term of $q_1$ and $q_2$ respectively. Now if we
expand the product $q_1q_2$ but do not collect terms, then
$a_{j_1}b_{j_2}\lambda^{j_1+j_2}$ is the unique term with the
least order. So the order of $c_{j_1+j_2}\lambda^{j_1+j_2}$ has to
equal the order of $a_{j_1}b_{j_2}\lambda^{j_1+j_2}$. This implies
that the $\lambda$-initial term of $q_1q_2$ is
$c_{j_1+j_2}\lambda^{j_1+j_2}$. The assertion for $\pt$ in the
lemma hence follows from the fact that $j_1+j_2=0\Leftrightarrow
j_1=0 \text{ and } j_2=0$. (Remember that $j_1,j_2\ge 0$).
\end{proof}
 A direct consequence of the above lemma is the following corollary.
 \begin{cor}\label{c-relprime}
If $1/q_1(\lambda)$ is $\PT^\rho$ in $\lambda$ and
$1/q_2(\lambda)$ is $\NT^\rho$ in $\lambda$, then $q_1(\lambda)$
and $q_2(\lambda)$ cannot have a nontrivial common divisor in
$K[\lambda]$, i.e., they are relatively prime.
\end{cor}

\begin{dfn}
If $q(\lambda)$ can be factored as $q_1(\lambda)q_2(\lambda)$ such
that $1/q_1(\lambda)$ is $\PT^\rho$ in $\lambda$ and
$1/q_2(\lambda)$ is $\NT^\rho$ in $\lambda$, then we say that
$q(\lambda)$ is $\rho$-{\em factorable}, and
$q(\lambda)=q_1(\lambda)q_2(\lambda)$ is a $\rho$-{\rm
factorization}. Such factorization is unique (if it exists) up to
a constant in $K$.
\end{dfn}

\begin{thm}
Let $p(\lambda), q(\lambda)\in K[\lambda]$. If $q(\lambda)$ is
$\rho$-factorable, then $\CT^\rho_\lambda p(\lambda)/q(\lambda)$
is in $K$, i.e., is rational.
\end{thm}
\begin{proof}
Suppose $q(\lambda)=q_1(\lambda)q_2(\lambda)$ is such a
$\rho$-factorization. Since  $1/q_1(\lambda)$ is $\PT^\rho $ in
$\lambda$ and $1/q_2(\lambda)$ is $\NT^\rho$ in $\lambda$,
$q_1(\lambda)$ and $q_2(\lambda)$ are relatively prime in
$K[\lambda]$ by Corollary \ref{c-relprime}. Thus we have the
unique partial fraction expansion in $K(\lambda)$:
\begin{align} \label{e-3-frac-r}
\frac{p(\lambda)}{q(\lambda)}
=p_0(\lambda)+\frac{p_1(\lambda)}{q_1(\lambda)}
+\frac{p_2(\lambda)}{q_2(\lambda)},
\end{align}
where $p_i$ are polynomials in $\lambda$ for $i=0,1,2$ and
 $\deg p_i(\lambda)< \deg q_i(\lambda)$ for
$i=1,2$. Since when expanded as series in $\CC^\rho\ll
\lambda,\mb{x}\gg$, $p_0(\lambda)$ and $p_1(\lambda)/q_1(\lambda)$
contains only nonnegative powers in $\lambda$, and
$p_2(\lambda)/q_2(\lambda)$ contains only negative powers in
$\lambda$, we have
$$\pt_\lambda \uprho \frac{p(\lambda)}{q(\lambda)} =
p_0(\lambda)+\frac{p_1(\lambda)}{q_1(\lambda)}.$$ Thus
$\ct^\rho_\lambda=p_0(0)+p_1(0)/q_1(0)$ is in $\CC(\mb{x})$.
\end{proof}

This theorem generalizes a result of Hadamard \cite[Proposition
4.2.5 ]{EC1}, which says that the Hadamard product of two rational
power series is rational. This statement can be easily seen from
the following observation: Let $f=\sum_{k\ge 0} f_k x^k$ and
$g=\sum_{k\ge 0} g_k x^k$. Then the Hadamard product of $f$ and
$g$ is
$$\sum_{k\ge 0}f_kg_kx^k = \ct_\lambda f(\lambda) g(x/\lambda),$$
where we are taking the constant term in $\CC\ll \lambda, x\gg$
for the RHS of the above equation.

For any  total ordering $\preceq^\rho$ on the monomials of
$K(\lambda)$, we let $\preceq^{\bar{\rho}}$ be the total ordering
such that $m_1\prec^{\bar{\rho}} m_2$ if and only if
$m_1\succ^\rho m_2$ for all monomials $m_1$ and $m_2$.

Then we have a sort of reciprocity invariant, for which we need
three more notations. We use the notation $\ct_{\lambda=0}
F(\lambda)$ to indicate that $F(\lambda)$ is treated as an element
in $K((\lambda))$ and $\ct_{\lambda=\infty} F(\lambda)$ to
indicate that $F(\lambda)$ is treated as an element in
$K((\lambda^{-1}))$. We define
\begin{align}\label{e-def-inv}
\mathcal{I}_\lambda F(\lambda) =\ct_{\lambda=0} F(\lambda)
+\ct_{\lambda=\infty} F(\lambda).
\end{align}

\begin{thm}\label{t-stan-inv}
Suppose that $p(\lambda),q(\lambda)\in K[\lambda]$, with
$q(\lambda)$ being $\rho$-factorable. Then the following is always
true as rational functions in $K$:
\begin{align}\label{e-stan-inv}
\ct_\lambda \uprho \frac{p(\lambda)}{q(\lambda)} +\ct_\lambda
\mbox{} ^{\bar{\rho}}\, \frac{p(\lambda)}{q(\lambda)}
=\mathcal{I}_\lambda \frac{p(\lambda)}{q(\lambda)}.
\end{align}
\end{thm}

Theorem \ref{t-stan-inv} gives an invariant of a rational function
when taking the constant term in $\lambda$. This invariant is
independent of the choice of the total ordering $\le^\rho$ when
applicable. This fact is the key in our new approach to the
monster reciprocity theorem in Section 4. 

\begin{proof}[Proof of Theorem \ref{t-stan-inv}]
Write $q(\lambda)$ as $q_1(\lambda)q_2(\lambda)\lambda^s$, such
that $1/q_1(\lambda)$ is $\PT^\rho$ in $\lambda$,
$1/q_2(\lambda)$ is $\NT^\rho$, and $q_1(0)q_2(0)\ne 0$. Clearly
$s\ge 0$ and the partial fraction decomposition of
$p(\lambda)/q(\lambda)$ can be written as
$$\frac{p(\lambda)}{q(\lambda)}= \frac{p_{-1}(\lambda)}{\lambda^s}
+p_0(\lambda)
+\frac{p_1(\lambda)}{q_1(\lambda)}+\frac{p_2(\lambda)}{q_2(\lambda)},
$$ where $\deg p_{-1}<s$, $\deg p_1<\deg q_1$, $\deg p_2<\deg q_2$, and
$p_0$ is a polynomial.

Now we are going to apply different operators on this partial
fraction decomposition. Applying $\ct_\lambda^\rho$ to
$p(\lambda)/q(\lambda)$ gives us $p_0(0)+p_1(0)/q_1(0)$, and
applying $\ct_\lambda^{\bar{\rho}}$ to $p(\lambda)/q(\lambda)$
gives us $p_0(0)+p_2(0)/q_2(0)$. Therefore
$$\ct_\lambda \uprho \frac{p(\lambda)}{q(\lambda)}
+\ct_\lambda \mbox{} ^{\bar{\rho}} \frac{p(\lambda)}{q(\lambda)} =
2p_0 (0)+\frac{p_1(0)}{q_1(0)}+\frac{p_2(0)}{q_2(0)}.
$$

Applying $\ct_{\lambda=0}$ to $p(\lambda)/q(\lambda)$ gives us
$p_0(0)+p_1(0)/q_1(0) +p_2(0)/q_2(0)$, and applying
$\ct_{\lambda=\infty}$ to $p(\lambda)/q(\lambda)$ gives us $p_0(0)$.
Thus the theorem follows.
\end{proof}

\begin{rem}
In the proof of Theorem \ref{t-stan-inv}, we see that $\bar{\rho}$
can be replaced with $\sigma$ if $\sigma$ switches the $\pt$ and
$\nt$ properties of $1/q_1(\lambda)$ and $1/q_2(\lambda)$ with
respect to $\rho$.
\end{rem}


As an element of $K[\lambda]$, $q(\lambda)$ can be factored into
the product of irreducible polynomials. Let
$q(\lambda)=q_1(\lambda)\cdots q_k(\lambda)$ be such a
factorization. By Lemma \ref{l-3-property} $q(\lambda)$ is
$\rho$-factorable if and only if every $1/q_i$ is either
$\PT^\rho$ or $\NT^\rho$. When this is the case, the
$\rho$-factorization can be obtained by collecting similar terms.

Elliott-rational functions are $\rho$-factorable for any $\rho$.
Such a function $F$ can be written as follows:
\begin{equation}\label{e-2-Mac-F2}
F=\frac{p(\lambda)}{(\lambda^{j_1}-a_1)\cdots (\lambda^{j_n}-a_n)
(\lambda^{k_1}-b_1)\cdots (\lambda^{k_m}-b_m)},
\end{equation}
where $p(\lambda)$ is a polynomial in $\lambda$, $j_i$ and $k_i$
are positive integers, $m$ and $n$ are nonnegative integers, and
$a_i$ and $b_l$ are monomials independent of $\lambda$. For a
particular $\rho$, we require that $1/(\lambda^{j_i}-a_i)$ is
$\NT^\rho$ in $\lambda$, and $1/(\lambda^{k_i}-b_i)$ is $\pt^\rho$
in $\lambda$. Note that $a_1$ can be $0$. ``The method of Elliott"
\cite[p.~111--114]{mac} shows that
 $\ct^\rho_\lambda F$ is always Elliott-rational.

A rational function of $\lambda$ is \emph{proper} in $\lambda$ if
the degree in $\lambda$ of its numerator is less than that of its
denominator.
\begin{cor}\label{c-macrec}
Let $F(\lambda)$ be of the form \eqref{e-2-Mac-F2}. If $F(0)=0$,
and $F(\lambda)$ is proper in $\lambda$, then for any $\rho$, we
have a reciprocity formula
$$\ct_\lambda \uprho F(\lambda) =-\ct_\lambda \mbox{}^{\bar{\rho}} \,
F(\lambda),$$ where both sides are regarded as elements in $K$.
\end{cor}

More generally, a rational function $F$ is said to have the
$R$-\emph{property} with respect to $\rho$ if
\begin{align}
\ct_\Lambda \mbox{}^{\rho} F =(-1)^d \ct_\Lambda
\mbox{}^{\bar{\rho}} F.
\end{align}
for some integer $d$. Here we restrict our interest to the case
when $d$ equals $r$, the number of $\lambda$'s. We have the
following reciprocity formula for Elliott-rational functions.

\begin{thm}\label{t-rec-erroterm} Let $F(\Lambda,\mb{x})$ be an
Elliott-rational function and let $\le ^\rho$ be a total ordering
on $\ZZ^{n+r}$ that is compatible with its additive group
structure. Then
\begin{align}
\ct_{\Lambda} \mbox{}^{\bar{\rho}} F =(-1)^{r} \ct_\Lambda \uprho
F + \sum_{i=0}^{r-1} (-1)^i {\ct_{\lambda_{r},\dots,
\lambda_{i+2}}}^{\!\!\!\!\!\!\rho} \;\mathcal{I}_{\lambda_{i+1}}
\ct_{\lambda_{i},\dots,
  \lambda_{1}}\mbox{}^{\!\!\!\!\bar{\rho}} \, F,
\end{align}
where $\ct_{\lambda_i,\dots, \lambda_1}^{\;\rho}$ is the identity
operator for $i=0$ and similar for $\ct_{\lambda_{r},\dots,
  \lambda_{i+2}}^{\;\bar{\rho}} $ when $i=r-1$.
\end{thm}
\begin{proof}
Since we are always taking constant terms in $\lambda_i$, we omit
the $\lambda$ for convenience. We compute the following in two
different ways.
\begin{align}\label{e-middle-term}
\sum_{i=0}^{r-1} (-1)^i \ct_{r,\dots,i+1} {}^{\!\!\!\!\bar{\rho}}
\; \ct_{i,\dots,
  1}{}^{\!\rho} F +(-1)^i  \ct_{r,\dots,
  i+2}{}^{\!\!\!\!\bar{\rho}} \;\ct_{i+1,\dots,1} {}^{\!\!\!\!\rho} F.
\end{align}

Using Theorem \ref{t-stan-inv}, we can rewrite
\eqref{e-middle-term} as
$$ \sum_{i=0}^{r-1} (-1)^i \ct_{r,\dots,
  i+2}{}^{\!\!\!\!\bar{\rho}}\; \mathcal{I}_{i+1}\; \ct_{i,\dots, 1}{}^{\!\rho} F.$$

On the other hand, most of the terms in \eqref{e-middle-term}
cancel with each other. The only terms left are given by
$$\ct_{r,\dots,2,1} {}^{\!\! \!\bar{\rho}} F +(-1)^{r-1}
\ct_{r,\dots,2,1} {}^{\!\!\!\rho} F .$$ The proposition then
follows.
\end{proof}

Theorem \ref{t-rec-erroterm} gives the error term of the
reciprocity formula. A different error term representation was
given in \cite{stanley-local} in terms of cohomology. However, the
computation of this error term saved only a little work for
general $r$. Our formula for the error term is true for any fixed
order of $\lambda_1,\dots,\lambda_r$, and any fixed order of
$x_1,\dots,x_n$. This suggests that some simplifications might
exist and a better formula is possible. We have not succeeded in
finding a such formula.

Since simple equivalent condition for $F$ to have the R-property
is unlikely, we search for a sufficient condition. Corollary
\ref{c-Iproperty} and Proposition \ref{p-Rproperty1} below play
important roles in our formulating the monster theorem.

A rational function $F$ is said to have the $I$-\emph{property}
with respect to $\rho$ if for $i=1,2,\dots, r$, we have
\begin{align}
\mathcal{I}_{\lambda_{i}} \ct_{\lambda_{i-1}}{}^\rho \cdots \ct_
  {\lambda_{1}}{}^{{\rho}} F =0.
\end{align}

\begin{cor}\label{c-Iproperty}
If an Elliott-rational function has the I-property, then it has
the R-property.
\end{cor}

This result is a direct consequence of Theorem
\ref{t-rec-erroterm}. The special case of this corollary when the
Elliott-rational function has a monomial numerator was shown by a
complicated computation
 in \cite[Lemma
9.2]{stanley-rec}.

In the case $r=1$, Theorem \ref{t-rec-erroterm} gives the
equivalence between the I-property and the R-property. Moreover,
we have, as shown below, a nice equivalent condition
\cite[Proposition 10.3]{stanley-rec} for the R-property that
contains no algebraic expression.

\begin{prop}[\cite{stanley-rec}]\label{p-Rproperty1}
Let $\mathcal{E}(\mb{x};b)$ be the crude generating function
associated to an LD-system consisting of a single equation
$A\alpha=a_1\alpha_1+\cdots +a_n\alpha_n=b$:
$$\mathcal{E}(\mathbf{x};b)=\frac{\lambda^{-b}}{(1-\lambda^{a_1}x_1)\cdots (1-\lambda^{a_n}x_n)}.$$
Then the following four conditions are equivalent for any $\rho$:
\begin{enumerate}
    \item $\mathcal{E}(\mb{x};b) $ has the R-property.
    \item $\mathcal{E}(\mb{x};b) $ has the I-property.
    \item $\ct_{\lambda=0}\mathcal{E}(\mb{x};b)=0 $ and $\ct_{\lambda=\infty}\mathcal{E}(\mb{x};b)=0
    $.
    \item The following two conditions are both satisfied:
    \begin{enumerate}
        \item There does not exist a $\beta\in \ZZ$ with $A\beta=b$ such
    that $\beta_e<0$ if $a_e>0$ and $\beta_e\ge 0$ if $a_e<0$.
        \item There does not exist a $\gamma \in \ZZ$ with $A\gamma =b$
    such that $\gamma_e\ge 0$ if $a_e>0$ and $\gamma_e<0$ if
    $a_e<0$.
    \end{enumerate}
    \end{enumerate}
\end{prop}
The proof of this proposition, which is not given in full detail
here, proceeds by showing that $\ct_{\lambda=0}
\mathcal{E}(\mb{x};b)$ and $\ct_{\lambda=\infty}
\mathcal{E}(\mb{x};b) $ have no common terms when expanded as
Laurent series. The reader is referred to \cite[Proposition
10.3]{stanley-rec} for details.

%
%
\section{The Monster Reciprocity Theorem}

Consider an LD-system $ A\alpha =\mb{b}$ as in \eqref{e-ax-b}. The
crude generating function $\mathcal{E}(\Lambda,\mb{x};\mb{b})$ is
an Elliott-rational function with a monomial numerator. We say
such a function has the \emph{matrix form} since we are going to
represent it by a matrix. The problem is to find a simple
sufficient condition for $\mathcal{E}(\Lambda,\mb{x};\mb{b})$ to
have the R-property. A homology version solution can be found in
\cite{stanley-local}. The best known result was Stanley's monster
reciprocity theorem \cite[Theorem 10.2]{stanley-rec}, which says
that the LD-system has the R-property if certain linear
combinations of its equations have the R-property. We present here
a simple approach to this problem.

The central idea of our approach to this problem, as in
\cite{stanley-rec}, is to apply Corollary \ref{c-Iproperty} and
Proposition \ref{p-Rproperty1}. If the following checking
procedure returns a true, then
$\mathcal{E}(\Lambda,\mb{x};\mb{b})$ has the I-property and hence
the R-property with respect to $\rho$. Note that the converse of
this statement is false.

The checking procedure for $\mathcal{E}(\Lambda,\mb{x};\mb{b})$:
\begin{enumerate}
    \item Let $T_1=\mathcal{E}(\Lambda,\mb{x};\mb{b})$. If $\mathcal{I}_{\lambda_1} T_1(\Lambda) \ne 0$ then \textsc{return} false.
    \item Write $\ct_{\lambda_1}\mbox{}^\rho  T_1(\Lambda)$ as
    a sum of matrix forms in an efficient way. For every matrix form $T_2$, if $\mathcal{I}_{\lambda_2}T_2\ne 0$,
    then \textsc{return} false.
    \item Repeat the above step for every matrix form $T_2$ with
    respect to $\lambda_2,$ and then for every $T_3$ with respect to $\lambda_3$,
   $\dots ,$ until we have checked if $I_{\lambda_{r}} T_r(\lambda_r)\ne 0$.
    If no false is returned, then \textsc{return} true.
\end{enumerate}
The basic tool in finding these $T_i$'s is partial fraction
decomposition of rational functions. Using the constant term
operators seems neater than using residue operators as in
\cite{stanley-rec}.

Our task is to find a simple equivalent condition for the checking
procedure to return a true. Such a condition will be our monster
reciprocity theorem. In order to do so, we represent a matrix form
$T$ as an augmented matrix. In fact, we can keep track everything
by adding a row of monomials in the $x$'s on the top and a column
of monomials in the $\lambda$'s to the left of an LD-system.
Therefore, the checking procedure will be done by matrix
operations. Note that using matrix operations is one important
aspect of the monster reciprocity theorem.

We use the following identification:
$$T= \frac{y_{n+1}\Lambda^{-\mb{b}}}{\prod_{i=1}^n (1-\Lambda^{C_i}y_i)}
\equiv \left[\begin{array}{ccc|c}
  y_1 & \cdots  & y_n & y_{n+1} \\
  C_1 & \cdots & C_n & \mb{b}
\end{array}\right], $$ where $y_i$ are monomials in $\mb{x}$, and $C_i$ are column
vectors. It would be clearer if we add $\lambda_i$ to the left of
the $i$th row, but this is unnecessary after applying Lemma
\ref{l-3-monomial} and requiring that the $i$th row (with $i\ge
2$) of $T_s$ is indexed by $\lambda_{s+i-1}$.

The row operations we are going to perform will never involve the
top row. The column operations, when acting on the first row, are
treated as multiplications instead of additions for the obvious
reason. We alow fractional entries and fractional powers. Roots of
unity might appear, but will not be a trouble.

Three special matrix operations will be useful. We define
$T\leftarrow C\langle i \rangle $ to be the matrix obtained from
$T$ by adding $-a_{1,j}/a_{1,i}$ times the $i$th column to the
$j$th column for all $j\ne i$. This operation is exactly Gaussian
column elimination by taking the $(2,i)$th entry of $T$ as the
pivot. Similarly we define the Gaussian row elimination
$T\leftarrow R \langle i \rangle $. The third operation
$T\leftarrow D\langle i\rangle $ is defined to be the matrix
obtained from $T$ by deleting the second row and the $i$-th
column.

Combination of the operations will also be used from left to
right. For instance, $T\leftarrow CR\langle i\rangle :=T\leftarrow
C\langle i \rangle\leftarrow R\langle i\rangle$. Since row
operations commute with column operations, we have $T\leftarrow
CR\langle i \rangle= T \leftarrow RC \langle i \rangle$. It is
easy to verify the following.
$$T\leftarrow CRD\langle i\rangle =T\leftarrow RCD\langle i\rangle
=T\leftarrow CD\langle i\rangle .
$$

For example, if $T$ is given by
\begin{align}\label{e-example}
T=\frac{\lambda_1^{-b}\lambda_2^{-c}}{(1-\lambda_1^3x_1/\lambda_2)(1-\lambda_2x_2/\lambda_1)
(1-x_3/\lambda_1^2\lambda_2)} \equiv \left[\begin{array}{ccc|c}
x_1 & x_2 & x_3 & 1 \\
3 & -1 & -2 & b \\
-1 & 1 & -1 & c
\end{array}\right],
\end{align}
 then
$$T\leftarrow C\langle 1\rangle  = \left[\begin{array}{ccc|c}
x_1 & x_2x_1^{\frac{1}{3}} & x_3 x_1^{\frac23}& x_1^{-\frac{b}{3}} \\
3 & 0 & 0 &0 \\
-1& \frac{2}3 & -\frac{5}{3} & c+\frac{b}{3}
\end{array}\right], \quad T\leftarrow R\langle 1\rangle   \left[\begin{array}{ccc|c}
x_1 & x_2 & x_3 & 1 \\
3 & -1 & -2 & b \\
0 & \frac{2}{3} & -\frac{5}{3} & c+\frac{b}{3}
\end{array}\right], $$
and
$$ T \leftarrow CD\langle 1\rangle  = T \leftarrow CRD\langle 1\rangle= \left[\begin{array}{cc|c}
 x_2x_1^{\frac{1}{3}} & x_3 x_1^{\frac23}& x_1^{-\frac{b}{3}} \\
 \frac{2}3 & -\frac{5}{3} & c+\frac{b}{3}
\end{array}\right].$$

The above three operations are generalized to sequences of
integers. For instance, $T\leftarrow R\langle i_1,\dots,i_p\rangle
$ is the matrix obtained from $T$ by applying Gaussian row
elimination by first taking the $(2,i_1)$th entry of $T$ as the
pivot, then taking the $(3,i_2)$th entry as the pivot, and so on.
However, the elimination cannot go backwards. For instance, we are
not allowed to eliminate the nonzero entries in the second row
when taking the $(3,i_2)$th entry as the pivot.

More precisely, pick out the $i_1,\dots,i_p$th columns of $T_1$,
and rearrange them as follows:
$$T_1(i_1,\dots,i_p):=\left[\begin{array}{cccc}
x_{i_1} & x_{i_2} & \cdots & x_{i_{p}} \\
a_{1,i_1} & a_{1,i_2} & \cdots & a_{1, i_{p}}\\
a_{2,i_1} & a_{2,i_2}  & \cdots & a_{2,i_{p}}\\
\vdots & \vdots & \vdots & \vdots\\
a_{r,i_1} & a_{r,i_2} & \cdots & a_{r, i_{p}} \\
\end{array}\right].$$
If all the pivots encountering are nonzero, then when ignoring the
top row
\begin{align}
\label{e-rowoperation} T_1(i_1,\dots,i_p)\leftarrow R\langle
i_1,\dots,i_p\rangle =\left[\begin{array}{cccc}
x_{i_1} & x_{i_2} & \cdots & x_{i_{p}} \\
a_{1,i_1} & a_{1,i_2} & \cdots & a_{1, i_{p}}\\
0 & a'_{2,i_2}  & \cdots & a'_{2,i_{p}}\\
\vdots & \vdots & \vdots & \vdots\\
0 & 0 & \cdots & a'_{p, i_{p}} \\
\mathbf{0} & \mathbf{0} & \cdots & \mathbf{0} \\
\end{array}\right]
\end{align}
will be an upper triangular square matrix followed by a zero
matrix, and \begin{align}
\label{e-TRC}T_1(i_1,\dots,i_p)\leftarrow RC\langle
i_1,\dots,i_p\rangle =\left[\begin{array}{cccc}
x_{i_1} & y_{i_2} & \cdots & y_{i_{p}} \\
a_{1,i_1} & 0 & \cdots & 0\\
0 & a'_{2,i_2}  & \cdots & 0\\
\vdots & \vdots & \vdots & \vdots\\
0 & 0 & \cdots & a'_{p, i_{p}} \\
\mathbf{0} & \mathbf{0} & \cdots & \mathbf{0} \\
\end{array}\right]
\end{align}
will be a diagonal square matrix followed by a zero matrix. Since
the matrix operations we have performed do not change the
determinants,
 the $a'_{s,i_s}$ can be inductively computed by the
formulas $a'_{1,i_1}=a_{1,i_1}$, and $\prod_{j=1}^s a'_{j,i_j}=\det
(a_{k,i_l})_{1\le k,l \le s}$.

We denote by $\mathcal{M}(z_1,\dots,z_k)$ is a generic monomial in
$z_1,\dots ,z_k$ whose exact expression is not needed.

Though we can formulate the monster reciprocity theorem for any
$\rho$, the result seems nicer if we assume that $\rho$ satisfies
the following condition:
$$ \forall y=\mathcal{M}(\mb{x}), \quad y\succ^\rho \lambda_s
\Rightarrow  y \mathcal{M}(\lambda_{s+1},\lambda_{s+2},\dots)
\succ\!\!\succ^\rho \lambda_s, \eqno{(\divideontimes)}.$$ For
example, condition $(\divideontimes)$ holds for any injective
$\rho$ such that $\rho(x_i)$ is a monomial in $x$ and
$\rho(\lambda_i)$ is a monomial in $\Lambda$. We will explain two
such $\rho$ in detail in the next section.

\begin{dfn}\label{def-contrib} With notation as in \eqref{e-TRC}, we define
$(i_1,\dots ,i_p)$ of distinct entries ranging from $1$ to $n$ to
be a \emph{contribution sequence} of length $p$ with respect to
$\rho$ if $y_{i_s}^{\text{sign}(-a'_{s,i_s})} \prec^\rho
\lambda_s$ for all $s$. The empty sequence is a contribution
sequence of length $0$.
\end{dfn}
The name contribution sequence is in correspondence with the
``pole sequence" in \cite{stanley-rec}. The condition in this
definition will be replaced with simple ones for two special
$\rho$ in the next section.

\begin{thm}[Generalized Monster Reciprocity Theorem]\label{t-monster}
Let $T$ be a matrix form corresponding to an $r$ by $n$ matrix of
full rank, and let $\preceq^\rho$ be a total ordering on  the
group of monomials in $\Lambda $ and $\mb{x}$ satisfying condition
$(\divideontimes)$. If for every contribution sequence $(i_1,\dots
,i_{p})$ of $T$ with $p<r$, the second row of $T\leftarrow
RD\langle i_1,\dots ,i_{p}\rangle $ has the R-property, then $T$
has the R-property with respect to $\rho$.
\end{thm}

\begin{proof}[Proof of Theorem \ref{t-monster}]
For given $T$ and $\preceq^\rho$. The checking procedure will
return a true if and only if every $T_p$ encountered has the
property that $\mathcal{I}_{\lambda_{p}} T_p=0$, which is the same
as the condition that the second row of $T_p$ has the R-property
by Proposition \ref{p-Rproperty1}.

We claim that $T_p$ must be \emph{similar} to the following form
for some contribution sequence $(i_1,\dots ,i_{p-1})$:
\begin{align}\label{e-Tp}
T_1 \leftarrow CD \langle i_1,\dots, i_{p-1}\rangle =T_1
\leftarrow RCD \langle i_1,\dots, i_{p-1}\rangle.
\end{align} The term similar will be explained later.
Assuming $T_p$ be given by \eqref{e-Tp}, we can complete the proof
of the theorem as follows. We observe that the $C$ operations
after the $R$ operations do not affect the $(p+1)$st row (and
below) of $T_1$. See \eqref{e-rowoperation}. Therefore the second
row of $T_p$ is the same as the second row of $T_1 \leftarrow RD
\langle i_1,\dots, i_{p-1}\rangle$.

We prove the claim by induction on $p$. The claim is trivial for
$p=1$. Now assume the claim is true for $p=s$ and we need to show
that the claim is true for $p=s+1$.

By choosing appropriate positive integer $N$ and letting
$\lambda=\lambda_s^{1/N}$, (note that $T_{s+1}$ will be independt
of the choice of $N$), we can assume that
$$T_s(\lambda_s)=T_s'(\lambda)=\frac{\lambda^{-b} \tilde{y}_{m+1}}{\prod_{i=1}^m (1-\lambda^{a_i}\tilde{y}_i)}=
\left[\begin{array}{ccc|c}
y_1 &\cdots & y_m & y_{m+1}\\
a_1 &\cdots & a_m & b \\
* & \cdots & * & * \end{array}\right],$$ where $a_k$ and $b$ are integers,
$\tilde{y}_k =y_k \mathcal{M}(\lambda_{s+1},\lambda_{s+2},\dots)$,
and the $*$'s are column of integers that we do not care. Dividing
the second row of $T_s'$ by $N$ will give us the the second row of
$T_s$. Since we have deleted $s-1$ columns, $m$ equals $n-s+1$.

We observe that $y_k=x_{k'} \mathcal{M}(x_{i_1},\dots
,x_{i_{s-1}})$ for $k\le m$, where $k'-k$ equals the number of
$j$'s such that $k'> i_j$. Therefore $y_1,\dots ,y_m$ are
independent of each other. It is now straightforward to check that
the partial fraction decomposition of $T'_s(\lambda)$ with respect
to $\lambda$ is:
$$T_s'(\lambda)=L(\lambda) +\sum_{k=1}^m \sum_{j=1}^{|a_k|}
\frac{-\zeta_{k,j}\tilde{y}_k^{-1/a_k}}{\lambda-\zeta_{k,j}\tilde{y}_k^{-1/a_k}}\cdot
\frac{(\zeta_{k,j}\tilde{y}_k^{-1/a_k})^{-b} \tilde{y}_{m+1}}{a_k
\prod_{i=1, i\ne k}^m
(1-(\zeta_{k,j}\tilde{y}_k^{-1/a_k})^{a_i}\tilde{y}_i)} ,$$ where
$L(\lambda)$ is a Laurent polynomial in $\lambda$, and
$\zeta_{k,j}$ ranges over all ${a_k}$th roots of unity.

By Proposition \ref{p-Rproperty1}, $\mathcal{I}_{\lambda_s} T_s
=0$ implies that $\ct_{\lambda_s=\infty}T_s=0$ and hence
$\ct_{\lambda_s}^\rho L(\lambda)=0$. Together with the fact that
for any $u$ independent of $\lambda$,
$$\ct_{\lambda_s}{}^\rho \frac{1}{\lambda-u}=\left\{%
\begin{array}{ll}
    0, &  \text{ if } u \succ^\rho \lambda, \\
    (-u)^{-1}, & \text{ if } u \prec^\rho \lambda, \\
\end{array}%
\right.     $$ we have
$$\ct_{\lambda_s}\mbox{}^\rho T_s =\sum_k \sum_{j=1}^{|a_k|}
\frac{(\zeta_{k,j}\tilde{y}_k^{-1/a_k})^{-b} \tilde{y}_{m+1}}{a_k
\prod_{i=1, i\ne k}^m
(1-(\zeta_{k,j}\tilde{y}_k^{-1/a_k})^{a_i}\tilde{y}_i)},$$ where
the sum ranges over all $k$ such that
$\zeta_{k,j}\tilde{y}_k^{-1/a_k}\prec^\rho
\lambda=\lambda_s^{1/N}$, which, by condition $(\divideontimes)$,
is equivalent to $y_i^{\text{sign}(-a_i)} \prec^\rho \lambda_s$.
For such $k$, we can check that $(i_1,\dots, i_{s-1},k')$ is a
contribution sequence.

Let \begin{align*} \mathcal{T}_{k}(x_{k'}^{-1/a_k})
=\frac{(\tilde{y}_k^{-1/a_k})^{-b} \tilde{y}_{m+1}}{ \prod_{i=1,
i\ne k}^m (1-(\tilde{y}_k^{-1/a_k})^{a_i}\tilde{y}_i)} &= T_s
\leftarrow CD \langle k \rangle,
\end{align*} where we
emphasize $\mathcal{T}_k$ as a function of $x_{k'}^{-1/a_k}$. By
delaying the deletion procedure, we can check that
$$\mathcal{T}_{k}(x_{k'}^{-1/a_k})= T_1
\leftarrow CD \langle i_1,\dots ,i_{s-1}, k' \rangle. $$ Then
$\ct_{\lambda_s} T_s$ is a sum of $T_{s+1}$'s, each have the form
$$\mathcal{T}_k^{(j)} = \frac{1}{a_k} \mathcal{T}_{k}( \zeta_{k,j}x_{k'}^{-1/a_k})  $$
for some $j$ with $1\le j \le |a_k|$ and $k$ with $(i_1,\dots
,i_{s-1}, k')$ being a contribution sequence.

We say that $\mathcal{T}_k^{(j)}$ is similar to $\mathcal{T}_k$
and it is clear that $\mathcal{T}_k^{(j)}$ has the I-property if
and only if $\mathcal{T}_k$ has. This completes the proof of the
claim.


%
%

\end{proof}

\begin{rem}\label{rem-monster}
If $T$ satisfies the condition in Theorem \ref{t-monster}, then
for $p\le r$, it follows from the proof that $\ct_{\lambda_1,\dots
,\lambda_p}^\rho T$ can be expressed as a sum of group terms
indexed by contribution sequences $(i_1,\dots ,i_p)$ of $T$, with
the corresponding group being a sum of terms similar to
$T_1\leftarrow \langle i_1,\dots ,i_p \rangle$. In particular,
$\ct_{\Lambda}^\rho$ can be expressed as at most a sum of
$n(n-1)\cdots (n-r+1)$ groups, since there are at most
$n(n-1)\cdots (n-r+1)$ contribution sequences. A fast way to
compute the sum for each group and an effective way to reduce the
number of contribution sequences will result in an efficient
algorithm for computing $E(\mb{x},\mb{b})$.
\end{rem}

\section{Examples and Applications}
To apply Theorem \ref{t-monster} to a particular $LD$-system, we
need to choose a working field $\CC^\rho \ll \Lambda, \mb{x}\gg$
to work with. The choice of $\rho$ is not unique, but we will
concentrate on two special cases that always work. One is the case
that $\rho $ is the identity; the other is equivalent to that in
\cite{stanley-rec}. In both cases, we can simplify the condition
in finding the contribution sequence.

Case 1: Let $\CC\ll \Lambda, \mb{x}\gg$ be the working field. The
condition $y_{i_s} ^{\text{sign} (-a'_{s,i_s})} \prec \lambda_s $
in Definition \ref{def-contrib} can be replaced with
$\text{sign}(a'_{s,i_s}k_l)
> 0$, where if we write $y_{i_s}=x_1^{k_1}\cdots x_n^{k_n}$, then
$l$ is the largest such that $k_l\ne 0$. In practice, we put the
sign of $k_l$ at the upper front of  $y_j$.

\begin{exa}Let $(\mathbf{E},(b,c))$ the following LD-system:
\begin{align*}
3\alpha_1 -\alpha_2-2\alpha_3 &=b, \\
-\alpha_1+\alpha_2-\alpha_3 &=c.
\end{align*}
\end{exa}
Then the crude generating function $T=T_1$ of this LD-system is
given by \eqref{e-example}. Using Maple, we find that
$(\mathbf{E},(b,c))$ has the R-property for all $(b,c)$ plotted by
$\bullet$, and has I-property for all $(b,c)$ plotted by $\circ$
in the following Figure \ref{f-rec1}, where we tested all $-12\le
b,c\le 12$. Thus the R-property does not implies the I-property.
\begin{figure}[h]\label{f-rec1}
 $$ \includegraphics[width=10cm]{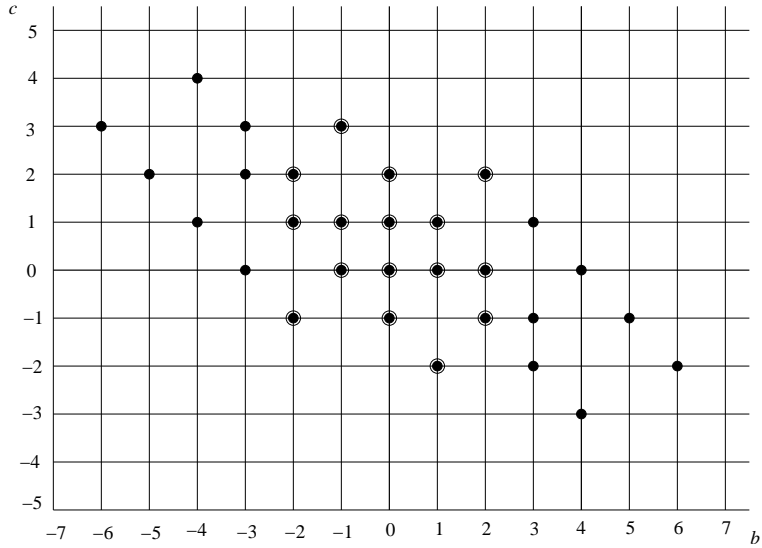}$$
  \caption{The R-property and I-property for $(\mathbf{E},(b,c))$.}\label{f-rec1}
\end{figure}

Let $\CC\ll \Lambda, \mb{x}\gg$ be the working field. We want to
apply Theorem \ref{t-monster} to find such pairs.

Since the second row of $T_1$ has only one positive entries, only
$(1)$ is a contribution sequence of length $1$. So after
eliminating $\lambda_1$, we get a sum of three terms similar to
$T_2$ given by
$$ T_2=T_1 \leftarrow CD\langle 1\rangle = \left[\begin{array}{ccc|c}
\oslash & ^+ x_2x_1^{\frac{1}{3}} & ^+ x_3 x_1^{\frac23}& x_1^{-\frac{b}{3}} \\
 \oslash& \frac{2}3 & -\frac{5}{3} & c+\frac{b}{3}
\end{array}\right],$$
where we kept the first column to keep track the original column
numbers. Now it is easy to see that the only contribution sequence
of length $2$ is $(1,2)$, though we do not need it.

Therefore, Theorem \ref{t-monster} tells us that the LD-system has
the R-property if the following two equations have the R-property:
 $$\begin{array}{rrrl}
 3\alpha_1 & -\alpha_2 & -2\alpha_3 &=b,\\
&2\alpha_2 & -5 \alpha_3 &=3c+b.
\end{array}$$ Using Maple, we find all such $(b,c)$ as plotted by $\circ$ in Figure \ref{f-rec1}:

\begin{exa}
Consider the equivalent LD-system $(\mb{E}'(-c,b))$:
\begin{align*}
\alpha_1-\alpha_2+\alpha_3 &=-c, \\
3\alpha_1 -\alpha_2-2\alpha_3 &=b,
\end{align*}
where we multiplied both sides of the second equation by $-1$ and
switched the two equations.
\end{exa} We need to find $(b,c)$ for $S$ to
have the R-property, where
$$S=\left[\begin{array}{ccc|c}
x_1 & x_2 & x_3 & 1 \\
1 & -1 & 1 & -c\\
3 & -1 & -2 & b
\end{array}\right].
$$
This time we have two contribution sequences of length $1$: $(1)$
and $(3)$. Therefore, Theorem \ref{t-monster} tells us that the
LD-system has the R-property if the following three equations have
the R-property, where the second and third equation are from
$S\leftarrow \langle 1\rangle$ and $S\leftarrow \langle 3\rangle$:
$$\begin{array}{rrrl}
\alpha_1 &-\alpha_2 &+\alpha_3 & =-c, \\
&2\alpha_2 &-5\alpha_3 & =b+3c, \\
5 \alpha_1& -3\alpha_2 & &=b-2c.
\end{array}
$$
Using Maple, we obtain the same pairs $(b,c)$ as in the previous
one, i.e., those plotted by $\circ$ in Figure \ref{f-rec1}. All
these three equations are needed to apply Theorem \ref{t-monster}.
The following coincidence is worth mentioning: if we only consider
the second and the third equation, we will get all $(b,c)$ plotted
by $\bullet$ in Figure \ref{f-rec1}, i.e., those $(b,c)$ such that
$(\mb{E}',(b,c))$ has the R-property. Since the first equation
comes from the empty contribution sequence, we come back to check
the previous example, which is obviously not the case.

Case 2: Let $\rho$ be the injective homomorphism into $\CC\ll
\mb{x}, \Lambda, t\gg$ by $\rho(x_i)=x_it$ and
$\rho(\lambda_i)=\lambda_{r-i+1}.$ Then the condition in
Definition \ref{def-contrib} can be replaced with $\text{sign}
(da'_{s,i_s})
>0$, where $d$ is the total degree of $y_{i_s}$ in the $x$'s.
Since we only need to keep track of the total degree of the $x$'s,
the $x_i$ in the top row of $T$ can be replaced with $1$. The
monster reciprocity theorem obtained this way is similar to that
of \cite{stanley-rec}, in which the computation used integration
along the circles $|\lambda_i|=1-\epsilon_i$ with $1>\!\!>
\epsilon_1>\!\!> \epsilon_2 >\!\!>\cdots$, where $>\!\!>$ means
``much greater", and the $x_i$ is taken to satisfy $|x_i|=\delta
<1$ for some positive real number $\delta $. In fact, the
condition as in Definition \ref{def-contrib} was completely
written in terms of determinants.

Detailed example for this case, which will not be given here, can
be found in \cite[p.~245]{stanley-rec}.

Now let us consider Linear homogeneous Diophantine system
(LHD-system for short). We shall use Theorem \ref{t-monster} to
derive the following theorem, which implies the reciprocal domain
theorem \cite[Proposition~8.3]{stanley-rec} including Theorem
\ref{t-3-recip}.

\begin{thm}\label{t-rec-dom}
Suppose $T$ is a matrix form corresponding to an LHD-system of
full rank. Then for any $\rho$ satisfying $(\divideontimes)$, $T$
has the R-property if and only if either $\ct_{\Lambda}^\rho
T=\ct_{\Lambda}^{\bar{\rho}} T=0$ or $\ct_{\Lambda}^\rho T\ne 0$
and $\ct_{\Lambda}^{\bar{\rho}} T\ne 0$.
\end{thm}
If we let $\rho$ be the identity map, then we get Theorem
\ref{t-3-recip}. If we let $\rho(x_i)=x_i$ for $i=1,\dots ,p$ and
$\rho(x_i)=x_i^{-1}$ for $i=p+1,\dots n$ for $p$ with $1<p<n$, and
$\rho(\lambda_i)=\lambda_i$ for all $i$, then we will get the
reciprocal domain theorem \cite[Proposition~8.3]{stanley-rec}.

\begin{proof}[Proof of Theorem \ref{t-rec-dom}]
If $T$ has the R-property, then
$$ \ct_{\Lambda}{}^\rho T=(-1)^r\ct_{\Lambda}{}^{\bar{\rho}} T.$$
The implication thus follows. Now we show the converse is true.

Obviously we can suppose $r>0$, $\ct_{\Lambda}^\rho T\ne 0$ and
$\ct_{\Lambda}^{\bar{\rho}} T\ne 0$. We first show that the second
row of $T$ has the R-property. Since $T$ corresponds to an
LHD-system, we can write
$$T= \frac{1}{\prod_{i=1}^n (1-\tilde{y}_i \lambda_1^{a_i} )},$$
where $\tilde{y}_i$ is a monomial independent of $\lambda_1$. If
some of the $a_i$ are positive and some of the $a_i$ are negative,
then Corollary \ref{c-macrec} applies and the second row of $T$
has the R-property. Otherwise, one of $\ct_{\lambda_1=0} T$ and
$\ct_{\lambda_1=\infty} T$ will be $0$ and the other will be
nonzero. (Note that since the LHD-system has full rank, the case
that $a_i=0$ for all $i$ will not happen.) The statement then
follows by Proposition \ref{p-Rproperty1}.

Now by Lemma \ref{l-3-monomial}, if $T'$ is obtained from $T$ by
elementary row operations, then $\ct_{\Lambda}^\rho T'\ne 0$ and
$\ct_{\Lambda}^{\bar{\rho}} T'\ne 0$. Therefore, the second row of
$T'$ has the R-property. This means every linear combination of
the equations of $T$ has the R-property. Thus the theorem follows
from Theorem \ref{t-monster}.
\end{proof}

\begin{rem}
The proof of the theorem, together with Remark \ref{rem-monster},
in fact shows the following statement: If $\ct_{\Lambda}^\rho T\ne
0$ and $\ct_{\Lambda}^{\bar{\rho}} T\ne 0$, then
$\ct_{\lambda_1,\dots,\lambda_p}^\rho T$ is proper in all
$\lambda_i$ for $i>p$. On the other hand, a simple proof of the
statement will lead to a simple proof of Theorem \ref{t-3-recip}.
If we restrict ourself in $\CC[\Lambda,\Lambda^{-1}][[\mb{x}]]$,
the best known proof of Theorem \ref{t-3-recip} should be that
given by the author in \cite{xinthesis}, which is included in the
next section.
\end{rem}

The above remark suggest a way to reduce the number of
contribution sequences of an LHD-system: Following the notation as
in Remark \ref{rem-monster}, since every $T_p$ has the R-property
for $\lambda_p$, we have a choice to choose all those terms with
contribution or all those terms (with a minus sign) without
contribution. The author is managing to develop a computer program
implementing these techniques.

\section{Linear Homogeneous Diophantine Systems}

We are concentrating on linear homogeneous Diophantine systems
(LHD-systems for short), i.e., $A\alpha=\mathbf{0}$. Recall that
$C_i$ is the $i$th column vector of $A$. We omit the $\mb{0}$ so
that $E$ and $\bar{E}$ are the sets of all solutions of
$A\alpha=\mathbf{0}$ in $\NN^n$ and $\PP^n$ respectively, and
similar for other notations. Since the proof closely related the
linear system and its associate generating functions, we restate
them as follows.

\begin{align}
\label{e-mac-Ebx1} E(\mb{x})=\ct_\Lambda \mathcal{E}(\mb{x}),
\qquad \qquad \bar{E}(\mb{x})=(-1)^n\ct_\Lambda
\mathcal{E}(\mb{x^{-1}}).
\end{align}

We are going to prove Proposition \ref{p-3-recip}, i.e., to show
that if $\bar{E}$ is nonempty, then
\begin{align}\label{e-3-recip1}
    \ct_\Lambda  \mathcal{E}(\mathbf{x}) =(-1)^{\mathrm{rank}(A)}\ct_\Lambda \uprho
    \mathcal{E}(\mathbf{x}),
    \end{align}
    where we are taking constant term of MN-series and $\rho(x_i)=x_i^{-1}$ for all $i$.

We shall see that all of the work is done algebraically. First,
let us see some facts. Exchanging column $i$ and $j$ corresponds
to exchanging $x_i$ and $x_j$. Row operations, which will not
change the solutions of $A\alpha=0$, are equivalent to multiplying
$A$ on the left by an invertible matrix. This fact can be obtained
by applying Lemma \ref{l-3-monomial}.

Let us see the simple case of
 $r=1$. In this case, $\mathcal{E}(\mb{x}) $ has the form:
$$ \mathcal{E}(\mb{x})=\prod_{i=1} ^n \frac{1}{1-\lambda^{a_i}x_i}.$$
The condition that $\bar{E} $ is nonempty is equivalent to saying
that some of $a_i$ have to be positive and some of $a_i$ have to
be negative. Thus when written in the normal form of a rational
function in $\lambda$, $\mathcal{E}(\mb{x})$ is proper and its
numerator divides $\lambda$. So Proposition \ref{p-3-recip}
follows from Corollary \ref{c-macrec}.

The general case does not seem to work along this line because of
two problems. One is how to use the conditions that $\bar{E}$ is
nonempty, and the other is how to connect to the rank of $A$. The
proof we are going to give uses induction and Elliott's reduction
identity \cite[p. 111--114]{mac}, which is easy to check and is
not given here.

Clearly if $a_{11},\dots ,a_{1,n}$ are all positive or are all
negative, then $\bar{E}$ is empty. So we can assume that
$a_{11}>0$ and $ a_{12}<0$. Applying Elliott's reduction identity
on $\lambda_1$, we get:
\begin{align*}
\mathcal{E}(\mb{x})&= \frac{1}{1-\Lambda^{C_1+C_2}x_1x_2}
\left(\frac{1}{1-\Lambda^{C_1}x_1}+\frac{1}{1-\Lambda^{C_2}x_2}-1
\right) \prod_{i\ge 3}\frac{1}{1-\Lambda^{C_i}x_i}
\end{align*}
Now expand $\mathcal{E}(\mb{x})$ according to the middle term, and
denote the resulting three summans by $\mathcal{E}_1$,
$\mathcal{E}_2$, and $\mathcal{E}_3$ respectively. We have
\begin{align}\label{e-3-E123}
\mathcal{E}(\mb{x})=\mathcal{E}_1(x_1,x_1x_2,x_3,\dots )+
\mathcal{E}_2(x_1x_2,x_2,x_3,\dots)-
\mathcal{E}_3(x_1x_2,x_3,\dots ).
\end{align}
Then these $\mathcal{E}_i$ are very similar to $\mathcal{E}$.
Correspondingly, they are associated to matrices, and hence
solution spaces that lie in $\NN^n$ and $\PP^n$. More precisely,
$\mathcal{E}_i$, $i=1,2,3$, are associated to $A_1=
(C_1,C_1+C_2,C_3,\dots,C_n)$, $A_2=(C_1+C_2,C_2,C_3,\dots,C_n)$,
and $A_3=(C_1+C_2,C_3,\dots,C_n)$ respectively. Thus $E_i,
E_i(\mb{x})$ and $\bar{E}_i,\bar{E}_i(\mb{x})$ are defined
correspondingly.

Now the matrix $A_1$ is obtained from $A$ by adding the second
column to the first; the matrix $A_2$ is obtained from $A$ by
adding the first column to the second. They are obtained from $A$
through a column operation. So the rank of $A_1$ and $A_2$ are
both equal to that of $A$. The rank of $A_3$ might not equal the
rank of $A$.

Applying $\ct_\Lambda$ and $(-1)^n\ct_\Lambda^{\rho}$ to
\eqref{e-3-E123} respectively, we get our key induction equations.
\begin{align}
E(\mb{x})&=E_1(x_1,x_1x_2,x_3,\dots )+ E_2(x_1x_2,x_2,x_3,\dots)-
E_3(x_1x_2,x_3,\dots ),\label{e-3-Ex123}\\
\bar{E}(\mb{x})&=\bar{E}_1(x_1,x_1x_2,x_3,\dots )+
\bar{E}_2(x_1x_2,x_2,x_3,\dots)\nonumber \\
&\qquad \qquad \qquad \qquad \qquad \quad +
(-1)^{\rank(A)-\rank(A_3)}\bar{E}_3(x_1x_2,x_3,\dots
).\label{e-3-Ebx123}
\end{align}

Looking more closely at these $E_i$, we can see that up to
isomorphism, $E_1$, $E_2$, and $E_3$ are obtained from $E$ by
intersecting the half spaces $\alpha_1\ge \alpha_2$, $\alpha_1\le
\alpha_2$, and the hyperplane $\alpha_1=\alpha_2$ respectively.
For instance, $(\alpha_1,\alpha_2,\dots,)$ belongs to $E$ with
$\alpha_1\ge \alpha_2$ if and only if
$(\alpha_1-\alpha_2,\alpha_2,\dots)$ belongs to $E_1$. Thus
Elliott's reduction identity in fact corresponds to a signed
decomposition of $E$. Equation \eqref{e-3-Ex123} and
\eqref{e-3-Ebx123} could be explained directly from geometry.

We need two more lemmas to give our proof of Proposition
\ref{p-3-recip}. We shall see that the condition on $\bar{E}$
plays an important role.

If $\bar{E}$ is nonempty, then $\dim E=\dim \bar{E}=n-\rank(A)$.
Clearly, the dimension of the solution space of $A\alpha=0$ is
$n-\rank(A)$. Let $\gamma\in \bar{E}$, and let $\Upsilon_1,\dots
,\Upsilon_{n-\rank(A)}$ be a $\ZZ$-basis of the solution space in
$\ZZ^n$ with $\Upsilon_1=\gamma$. Then for sufficiently large $m$,
$m \gamma +\Upsilon_1,\dots ,m\gamma+\Upsilon_{n-\rank(A)}$ will
be a linearly independent set in $\bar{E}$.

\begin{lem}\label{l-3-empty-23}
Suppose that $\bar{E}$ is nonempty, and that $\bar{E}_i$ is
defined as above for $i=1,2,3$. Then any two of the $\bar{E}_i$
being nonempty implies that they are all nonempty.
\end{lem}
\begin{proof}
Suppose that $\bar{E}_1$ and $\bar{E}_2$ are nonempty. Then we
have elements $\beta$ and $\gamma$ in $\bar{E}$ such that
$\beta=(\beta_1,\beta_2,\dots)$ with $\beta_1>\beta_2$ and
$\gamma=(\gamma_1,\gamma_2,\dots)$ with $\gamma_1<\gamma_2$. Then
$(\gamma_2-\gamma_1)\beta+(\beta_1-\beta_2)\gamma$ is in $\bar{E}$
with the first two entries being equal. This means $\bar{E}_3$ is
nonempty.

Suppose that $\bar{E_1}$ and $\bar{E}_3$ are nonempty. Then we
have elements $\beta$ and $\delta$ in $\bar{E}$ such that
$\beta=(\beta_1,\beta_2,\dots)$ with $\beta_1>\beta_2$ and
$\delta=(\delta_1,\delta_2,\dots)$ with $\delta_1=\delta_2$. Then
for sufficiently large $m$, $m\delta -\beta$ is in $\bar{E}$ with
the first entry being smaller than the second. This means
$\bar{E}_2$ is nonempty.

The case that $\bar{E_2}$ and $\bar{E}_3$ are nonempty is similar
to the previous case.
\end{proof}

\begin{lem}\label{l-3-rankall}
If all of the $\bar{E_i}$ are nonempty, then $\rank (A_3)=\rank
(A)$.
\end{lem}
\begin{proof}
By hypothesis, it is clear that $E$ is not contained in the
hyperplane $\alpha_1=\alpha_2$. Thus the intersection of $E$ with
the hyperplane has dimension $\dim E -1$. So $\dim E_3$ is also
$\dim E-1$ and the rank of $A_3$ equals $n-1-\dim E_3=\rank (A)$.
\end{proof}

\begin{proof}[Proof of Proposition \ref{p-3-recip}]
The base case, when $A$ is the zero matrix, is trivial.

By exchanging rows, we can assume that not all of the entries in
the first row of $A$ are zero. Moreover, since the entries can not
be all positive or negative, we can assume the first entry is
positive and the second is negative by exchanging columns.

We use induction on $S_1(A)$, which is defined to be the sum of
the absolute values of all the entries in the first row. Now the
above argument applies, and it is easy to see that
$S_1(A_i)<S_1(A)$ for $i=1,2,3$. Applying Lemma
\ref{l-3-empty-23}, we can reduce the seven cases of $E_i$ being
nonempty or not into the following four cases:

Case 1: only $\bar{E}_1$ is nonempty. Let $\beta$ in $\bar{E}$ be
such that $\beta_1>\beta_2$. We claim that all $\alpha$ with
$A\alpha=0$ satisfy the condition $\alpha_1
>\alpha_2$, so that $E_2(x_1x_2,x_2,x_3,\dots)$ equals $E_3(x_1x_2,x_3,\dots)$, and hence by induction we have
\begin{align*}
E(\mb{x})&=E_1(x_1,x_1x_2,x_3, \dots)\\
&= (-1)^{\rank(n-A_1)} \bar{E}_1(\xx _1,\xx _1 \xx_2, \xx
_3,\dots)=(-1)^{n-\rank(A)}\bar{E}(\mb{x^{-1}}).
\end{align*}

If the claim does not hold, then  $\alpha_1\le \alpha_2$. But for
sufficiently large $m$, $m\beta -\alpha$ will produce an element
in $\bar{E}_2$  or $\bar{E}_3$, a contradiction.

Case 2: only $\bar{E}_2$ is nonempty. This is similar to case 1.

Case 3: only $\bar{E}_3$ is nonempty. This means that $E$ is
contained in the hyperplane $\alpha_1=\alpha_2$. Thus
$$E_1(x_1,x_1x_2,x_3,
\dots)=E_2(x_1x_2,x_2,x_3, \dots)=E_3(x_1x_2,x_3,\dots),$$ and we
have
$$\rank(A_3)=n-1-\dim(E_3)=n-\dim(E)-1=\rank(A)-1.$$
So
\begin{align*}
E(\mb{x})&=E_3(x_1x_2,x_3,\dots)\\
&=(-1)^{n-1-\rank(A_3)}\bar{E_3}(\xx _1\xx _2,\xx
_3,\dots)=(-1)^{n-\rank(A)}\bar{E}(\mb{\xx}).
\end{align*}

Case 4: all of $\bar{E}_i$ are nonempty. By induction, we see that
$$E_i(\mb{x})=(-1)^{n-\rank(A_i)} \bar{E}_i(\mb{x^{-1}})$$
for $i=1,2$, and that
$$E_3(x_2,x_3,\dots)=(-1)^{n-1-\rank{A_3}}\bar{E}(\xx
_2,\xx_3,\dots).$$ From Lemma \ref{l-3-rankall},
$\rank(A_3)=\rank(A)$. Thus together with our key induction
equations \eqref{e-3-Ex123} and \eqref{e-3-Ebx123}, we get
\begin{align*}E(\mb{x})=&E_1(x_1,x_1x_2,x_3,\dots)
+ E_2(x_1x_2,x_2,x_3, \dots)-E_3(x_1x_2,x_3,\dots)\\
=&(-1)^{n-\rank(A)}\left(\bar{E}_1(\xx _1, \xx _1\xx _2,\xx
_3,\dots)\right.\\
&\qquad \qquad\qquad\qquad\qquad\left. +\bar{E}_2(\xx _1\xx_2, \xx
_2,\xx _3,\dots)+\bar{E_3}(\xx _1\xx
_2,\xx _3,\dots)\right)\\
=&(-1)^{n-\rank{A}}\bar{E}(\mb{x}).
\end{align*}
\end{proof}

{\bf Acknowledgement:} I am very grateful to Richard Stanley for
introducing me to his inspiring work.

\bibliographystyle{amsplain}

\providecommand{\bysame}{\leavevmode\hbox
to3em{\hrulefill}\thinspace}
\providecommand{\MR}{\relax\ifhmode\unskip\space\fi MR }
\providecommand{\MRhref}[2]{%
  \href{http://www.ams.org/mathscinet-getitem?mr=#1}{#2}
} \providecommand{\href}[2]{#2}

\end{document}